\documentstyle[12pt]{article}
\begin{document}
\title{ The simple method of distinguishing the underlying differentiable
structures of algebraic surfaces }
\author{Andrej Tyurin }
\date { }
\maketitle{}
\begin{center}
Steklov Institute of Mathematics \\ Moscow
\end{center}
% bstab
\def\iff{\Leftrightarrow}
\def\then{\Rightarrow}
\def\der{\frac{d}{dt}}
\def\derz{\frac{d}{dt}|_{t=0}}
\def\pf{{\bf Proof:}}
\def\lm{{\bf Lemma:}}
\def\da{\downarrow}
\def\beg{\begin{equation}}
\def\en{\end{equation}}
\def\beq{\begin{eqnarray}}
\def\enq{\end{eqnarray}}

\def\bbbr{{\rm I\!R}} %reelle Zahlen
\def\bbbn{{\rm I\!N}} %natuerliche Zahlen
\def\bbbm{{\rm I\!M}}
\def\bbbh{{\rm I\!H}}
\def\bbbk{{\rm I\!K}}
\def\bbbp{{\rm I\!P}}
\def\bbbl{{\rm I\!L}}
\def\bbbone{{\mathchoice {\rm 1\mskip-4mu l} {\rm 1\mskip-4mu l}
{\rm 1\mskip-4.5mu l} {\rm 1\mskip-5mu l}}}

\def\bbbc{{\mathchoice {\setbox0=\hbox{$\displaystyle\rm C$}\hbox{\hbox
to0pt{\kern0.4\wd0\vrule height0.9\ht0\hss}\box0}}
{\setbox0=\hbox{$\textstyle\rm C$}\hbox{\hbox
to0pt{\kern0.4\wd0\vrule height0.9\ht0\hss}\box0}}
{\setbox0=\hbox{$\scriptstyle\rm C$}\hbox{\hbox
to0pt{\kern0.4\wd0\vrule height0.9\ht0\hss}\box0}}
{\setbox0=\hbox{$\scriptscriptstyle\rm C$}\hbox{\hbox
to0pt{\kern0.4\wd0\vrule height0.9\ht0\hss}\box0}}}}

\def\bbbg{{\mathchoice {\setbox0=\hbox{$\displaystyle\rm G$}\hbox{\hbox
to0pt{\kern0.4\wd0\vrule height0.9\ht0\hss}\box0}}
{\setbox0=\hbox{$\textstyle\rm G$}\hbox{\hbox
to0pt{\kern0.4\wd0\vrule height0.9\ht0\hss}\box0}}
{\setbox0=\hbox{$\scriptstyle\rm G$}\hbox{\hbox
to0pt{\kern0.4\wd0\vrule height0.9\ht0\hss}\box0}}
{\setbox0=\hbox{$\scriptscriptstyle\rm G$}\hbox{\hbox
to0pt{\kern0.4\wd0\vrule height0.9\ht0\hss}\box0}}}}

\def\bbbq{{\mathchoice {\setbox0=\hbox{$\displaystyle\rm Q$}\hbox{\raise
0.15\ht0\hbox to0pt{\kern0.4\wd0\vrule height0.8\ht0\hss}\box0}}
{\setbox0=\hbox{$\textstyle\rm Q$}\hbox{\raise
0.15\ht0\hbox to0pt{\kern0.4\wd0\vrule height0.8\ht0\hss}\box0}}
{\setbox0=\hbox{$\scriptstyle\rm Q$}\hbox{\raise
0.15\ht0\hbox to0pt{\kern0.4\wd0\vrule height0.7\ht0\hss}\box0}}
{\setbox0=\hbox{$\scriptscriptstyle\rm Q$}\hbox{\raise
0.15\ht0\hbox to0pt{\kern0.4\wd0\vrule height0.7\ht0\hss}\box0}}}}
\def\bbbt{{\mathchoice {\setbox0=\hbox{$\displaystyle\rm
T$}\hbox{\hbox to0pt{\kern0.3\wd0\vrule height0.9\ht0\hss}\box0}}
{\setbox0=\hbox{$\textstyle\rm T$}\hbox{\hbox
to0pt{\kern0.3\wd0\vrule height0.9\ht0\hss}\box0}}
{\setbox0=\hbox{$\scriptstyle\rm T$}\hbox{\hbox
to0pt{\kern0.3\wd0\vrule height0.9\ht0\hss}\box0}}
{\setbox0=\hbox{$\scriptscriptstyle\rm T$}\hbox{\hbox
to0pt{\kern0.3\wd0\vrule height0.9\ht0\hss}\box0}}}}
\def\bbbs{{\mathchoice
{\setbox0=\hbox{$\displaystyle     \rm S$}\hbox{\raise0.5\ht0\hbox
to0pt{\kern0.35\wd0\vrule height0.45\ht0\hss}\hbox
to0pt{\kern0.55\wd0\vrule height0.5\ht0\hss}\box0}}
{\setbox0=\hbox{$\textstyle        \rm S$}\hbox{\raise0.5\ht0\hbox
to0pt{\kern0.35\wd0\vrule height0.45\ht0\hss}\hbox
to0pt{\kern0.55\wd0\vrule height0.5\ht0\hss}\box0}}
{\setbox0=\hbox{$\scriptstyle      \rm S$}\hbox{\raise0.5\ht0\hbox
to0pt{\kern0.35\wd0\vrule height0.45\ht0\hss}\raise0.05\ht0\hbox
to0pt{\kern0.5\wd0\vrule height0.45\ht0\hss}\box0}}
{\setbox0=\hbox{$\scriptscriptstyle\rm S$}\hbox{\raise0.5\ht0\hbox
to0pt{\kern0.4\wd0\vrule height0.45\ht0\hss}\raise0.05\ht0\hbox
to0pt{\kern0.55\wd0\vrule height0.45\ht0\hss}\box0}}}}
\def\bbbz{{\mathchoice {\hbox{$\sf\textstyle Z\kern-0.4em Z$}}
{\hbox{$\sf\textstyle Z\kern-0.4em Z$}}
{\hbox{$\sf\scriptstyle Z\kern-0.3em Z$}}
{\hbox{$\sf\scriptscriptstyle Z\kern-0.2em Z$}}}}

\def\gaaa{\ifmmode
              {{\mbox{\deu a}}}%
          \else${{\mbox{\deu a\ }}}$%
           \fi}
\def\gbbb{\ifmmode
              {{\mbox{\deu b}}}%
          \else${{\mbox{\deu b\ }}}$%
           \fi}
\def\ccc{\ifmmode
        {\bbbc}%
          \else${\bbbc\ }$%
          \fi}
\def\fff{\ifmmode
        {{\mbox{\bf F}}}%
          \else${{\mbox{\bf F\ }}}$%
          \fi}
\def\fxfx{\ifmmode
              {{\bbbf}_p[x] }%
          \else${{\bbbf}_p[x]\ }$%
          \fi}
\def\fffp{\ifmmode
              {{\mbox{\bf F}}_p }%
          \else${{{\mbox{\bf F}}}_p\ }$%
          \fi}
\def\fffq{\ifmmode
              {{\bbbf}_q }%
          \else${{\bbbf}_q\ }$%
          \fi}
\def\gmmm{\ifmmode
              {{\mbox{\deu m}}}%
          \else${{\mbox{\deu m\ }}}$%
           \fi}
\def\nnn{\ifmmode
        {\bbbn}%
          \else${\bbbn\ }$%
          \fi}

\def\nnno{\ifmmode
        {{\bbbn}_0}%
          \else${{\bbbn}_0\ }$%
          \fi}

\def\ooo{\ifmmode
              {{\mbox{\deu o}}}%
          \else${{\mbox{\deu o\ }}}$%
           \fi}
\def\gppp{\ifmmode
              {{\mbox{\deu p}}}%
          \else${{\mbox{\deu p\ }}}$%
           \fi}

\def\ppp{\ifmmode
              {{\mbox{\deu p}}}%
          \else${{\mbox{\deu p\ }}}$%
           \fi}
\def\cppp{\ifmmode
              {\cal P}%
          \else${\cal P\ }$%
           \fi}

\def\qqq{\ifmmode
              {\bbbq}%
          \else${\bbbq\ }$%
           \fi}
\def\zzz{\ifmmode
              {\bbbz}%
          \else${\bbbz\ }$%
          \fi}
\def\grrr{\ifmmode
              {{\mbox{\deu r}}}%
          \else${{\mbox{\deu r\ }}}$%
           \fi}
\def\rrr{\ifmmode
              {\bbbr}%
          \else${\bbbr\ }$%
          \fi}
\def\st{\ifmmode
              {\star\star\star}%
          \else${\star\star\star}$%
          \fi}

\newcommand{\nat}{\nnn}
\newcommand{\rat}{\qqq}
\newcommand{\gan}{\zzz}

%\newfont{\deu}{eufm10}
%\newcommand{\frak}{\deu}
\newcommand{\goth}{\frak}
\newcommand{\frak}{\bf}
\newcommand{\deu}{\bf}
\newcommand{\notag}{\nonumber}
\newcommand{\text}{\mbox}
\newcommand{\phant}{\hspace*{6em}}

%end bstab

\addtocounter{section}{-1}
  \[  \]
\section{ Introduction}

   The purpose of this preprint is to construct a new invariant of the smooth
structure of a simply connected 4-manifold M so called Spin-polynomials

 \[ \gamma^{g,C}_1 (2,c_1,c_2) \in S^{d_1} H^2(M, \bbbz)  {~~~~~}   (0.1) \]
and to show how to use it to compare the smooth structures of rational surfaces
and surfaces of general type.

This Spin-polynomial (0.1)  is the analogue of the original Donaldson
polynomial
$$ \gamma^g(2,c_1,c_2)  \in S^d H^2(M, \bbbz)  {~~~~~}  (0.1')  $$

and depends on one extra index C given by a so called  $ Spin^{ \bbbc}
$-structure on M . It specifies a lift of the Stiefel-Whitney class  $w_2(M)
\in H^2(M, \bbbz _2) $ to some integer class  $C \in H^2(M, \bbbz ) $ .

The technical basis of the construction of the Spin-polynomial is the same as
of the ordinary polynomials (0.1') .It is our aim to compare the properties of
polynomials (0.1) and (0.1')

Here our aim is to consider only the simplest version of invariants of such
type and to use in applications the simplest arguments of proofs (now standart
for Donaldson's stuff ).

Much more sophisticated constructions and a vague discussion of the properties
of these invariants are contained in forthcoming article [T 3].

A special case of our  construction has been used in the article [P-T], where
some important basic theorems were proved. By this reason we will follow the
english translation of this article in terminology and notation (see [P-T]).

As applications of our techniques we will prove the non existence of algebraic
fake planes,Hirzebruch surfaces or quadrics.(Here "fake" means "diffeomorphic
to   ...,but algebraically non equivalent").

I should like to thanks the SONDERFORSCHUNGSBEREICH 170 in G\"ottingen for its
kind hospitality, support and the use of its facilities.

   \section {$ Spin ^ \bbbc $-structure }

 The Stiefel-Whitney class $ w_2 (M) \in H^2(M, \bbbz_2 ) $ of a smooth,compact
4-manifold M is the characteristic class of the lattice $ H^2(M, \bbbz) $ with
its intersection form $ q_M $ . This means that for every $ \sigma \in H^2(M,
\bbbz) $ \[ \sigma^2 \equiv \sigma.w_2 (M){~} mod {~}2 {~~~~~~~} (1.1) \]

{\bf Definition 1.1. }  Let M be a smooth, simply connected,compact 4-manifold
,then a class $ C \in H^2(M, \bbbz) $ such that $ C \equiv w_2(M)   mod  2 $ is
called a $Spin^ \bbbc $-structure of M.Thus to equip M with a $Spin^ \bbbc
$-structure is the same as to lift up  $w_2 $ to an integer class.

The set of all $Spin^ \bbbc $-structures on M is the affine sublattice \[
H_w(M) = \{ \sigma \in H^2(M, \bbbz) \mid  \sigma \equiv w_2(M) {~}mod {~}2 \}
{~~~~}(1.2) \]

 {\bf Remarks} 1) For every  $C \in H_w(M) $ \[ C^2 \equiv b^+_2- b^- _2 = I
{~} mod {~} 8 {~~~~~~~~} (1.3) \]
  (here I is the index or signature of M)

2) The diffeomorphism group Diff M of M acts on  $H_w(M)$  by the affine
transformations.

3) the lattice  $ H^2(M, \bbbz) $ 1-transitive acts on  $H_w(M) $ by the
formula: for  $ C \in H_w(M) $ \[ \sigma(C)=C+2 \sigma {~~~~~} (1.4) \]
Hence, the choice of a $ C_0 $ of $  H_w(M)$ gives an identification :
$H_w(M)=H^2(M).$

For the rest of this paper the symbol $ L_ \sigma  $  for  $\sigma \in H^2(M,
\bbbz) $ will denote a complex line bundle with the first Chern class \[ c_1(
L_ \sigma )= \sigma {~~~~~~} (1.5)\]

If M is equipped with a Riemannian metric g the $Spin ^\bbbc $ -structure C on
M defines a pair of rank 2 Hermitian vector bundles $ W^+ $ and $ W^- $ such
that the complexification of the tangent bundle can be decomposed as a tensor
product

 \[ TM_ \bbbc=(W ^-)^*  \otimes W^+ {~~~~~} (1.6) \]

 with

\[  \Lambda^2 W^ {+,-} = L_C {~~~~~~~}(1.7) \]

Moreover,if we equip  $ L_C $ with any Hermitian connection $
\bigtriangledown_0 ,$ then for any U(2)-bundle E on M and any Hermitian
connection \[ a \in  {A  \cal } _h  (E) {~~~~~~} (1.8) \]
on E we have a coupled Dirac operator \[ D ^{C, \bigtriangledown _0 }_a :
\Gamma ^ {\infty }(E \otimes W^+ ) \rightarrow \Gamma^  { \infty }(E \otimes
W^-) {~~~~~~~~~} (1. 9) \]
and an integer number
\[ \chi _C (E) = ind D ^{C,\bigtriangledown _0 } _a {~~~~~~~} (1.10 ) \] which
is the index of the Fredholm operator (1.9).

  It is easy to see that this integer number doesn't depend on the continuous
parameters g, $ \bigtriangledown _0 $ and $ a $ and depends on the $ Spin ^
\bbbc $-structure $ C \in H _w (M)$  only.

If we change the $ Spin ^ \bbbc $ -structure from C to C' (we can consider it
as   a result of the action (1.4) of $ \delta =(C'- C) / 2 $ on C ) we have \[
\chi _ {C+ 2 \delta }  (E) = \chi _C (E) + c_1 . \delta -   \delta ( \delta + C
) {~~~~~~} (1.11) \] because of  \[ \chi _{C + 2 \delta } (E) = \chi _C ( E
\otimes L _ \delta ) {~~~~~~~} (1.12) \] and by the Atiyah - Singer formula \[
\chi _C (E) = c_1 (c_1 + C)/2 + 2 \chi _C (L_o) -c_2  {~~~~~~~~} (1.13) \]
.Here $ L_0 $ denotes the  trivial    line   bundle   on M   and $ c_1, c_2 $
are  the Chern- classes  of  E.

Finally   we shoud point out that special structures on M can define a
canonical $Spin ^ \bbbc $-structure: for example if M is the underlying
structure of a complex surface S, then there is the natural $Spin ^ \bbbc
$-structure $ C = c_1  (S) = - K_S $  given by the anticanonical class or if M
admits a symplectic structure  $ \omega $,then there exists the canonical set
of complex structures on the tangent bundle TM=E with the same first Chern
class  $c_1 (M, \omega ) $ which
is the natural $Spin ^ \bbbc $-structure for $(M, \omega )$.

 \section { The definition of Spin-polynomials }

  If a $Spin ^ \bbbc $ -4-manifold (M,C) is equipped with a Riemannian metric g
then   for every $U(2)$-- bundle E the gauge - orbit space
\[{ \cal B } (E) = {\cal A^*}_h (E) /{ \cal G}  \]
of irreducible connections contains the subspace
\[ {\cal M}^g (E) \subset {\cal B}(E) {~~~~~~~} (2.1) \]
 of anti self dual connections with respect to the Riemannian metric g.

  We can consider the subspace of $\cal M$ $^g (E)$:
\[ {\cal M}^{g,C} _1 (E) = \{ (a) \in {\cal M}^g (E) | rk {~} ker D ^
{C,\bigtriangledown _0} _a \geq 1 \} {~~~~~~} (2.2) \]

Analogously,
\[ {\cal M}^{g,C}_2 (E) = \{ (a) \in {\cal M}^g (E) | rk {~} ker D ^
{C,\bigtriangledown _0} _a \geq 2 \} {~~~~~~}  \]
and so on.

 If $ (a) \in {\cal M}^{g,C}_1 - {\cal M}^{g,C}_2 $ and the family of Dirac
operators is in "general position" near (a),  then the
fibre of the normal bundle to  ${\cal M}^{g,C}_1 (E)$ at (a) is given by
\[ (N_{{\cal M} _1 \subset {\cal M }}  ) _{{~} (a)} = Hom  ( ker D_a, coker D_a
)  {~~~} (2.3) \]
with  $ker D_a = \bbbc $,  $coker D_a = \bbbc^{1 - \chi _C (E)} $ (if the index
of the Dirac operators is not positive ).

 Thus the virtual (expected) codimension of $ {\cal M}^{g,C}_1 (E) $

\[ v.codim {\cal M}^{g,C} _1 (E) = 2 - 2 \chi _C (E) {~~~~~~~~~~~~~} (2.4) \]

 On the analogy of the Freed - Uhlenbeck theorem,which says that for generic
metric g the moduli space$ \cal M$$^g (E) $ (2.1) is a smooth manifold of the
expected dimension with regular ends (see Theorem 3.13 of [F-U])  the following
fact was proved  in section 3 of Ch.2 of [P-T] .
\[  \]

  {\bf Transversality Theorem}.For generic pair (g, $ \bigtriangledown _0) \in
{\cal S} \times \Omega^1 $  of metric and connection on $L_C$, the moduli space
${\cal M}^{g,C}_1 (E)$ is smooth outside $\cal M$$^{g,C} _2 (E)$ of expected
codimension (2.4).
\[  \]

Moreover, $\cal M$$^g (E)$ admits a natural orientation (see [D 1] and [K]).But
 $\cal M$$^{g,C} _1 (E) $ admits the special orientation becouse its normal
bundle  (2.4)  has a natural complex structure. This orientation is described
in details  in section 5 of Ch.2 of [P-T].

  Now, we need the usual restrictions on the topology of M.We will suppose

 \[  b^+ _2  (M)= 2 p_g (M) +1   \ \ \ \ (2.5) \]
to be odd.Then both v.dim $\cal M$$^g (E) = 2 d $  and v.dim $\cal M$$^ {g,C}
_1(E) = 2 d_1 $ must be even .

 To compute the value of  $ \gamma ^{g,C}_1(E) $ evaluated at an argument $
(\sigma_1,...,\sigma_{d_1}) $ we need to consider Donaldson's realisation of
this collection of 2-cycles as a collection of smoothly embedded Riemannian
surfaces
$ ( \Sigma_1,...,\Sigma_ {d_1}) $ which are in general position in the
following sense :

  1) Any two surfaces  $ \Sigma _i $ and $ \Sigma _j $ meet transversally. And
let \[ \{ m_1,...,m_N \}  {~~~~}  (2.6) \]  be the set of all points wich are
intersection points for some i and j.

  2) exactly two surfaces pass through any point of intersection $m_i $ so that
in the flag diagram

\[  \begin{array}{rcl}
   & \{ m_ i \in \Sigma _j \} & {~~~~~~~~}  \\
 & &  {~~~~~~~~~~~~~~~~~~~} (2.7) \\
 \swarrow & & \searrow   \\
\{ m_i,...,m_ N \}   & & {~} \{ \Sigma _1,...,\Sigma _{d_1}  \}
\end {array} \]
the projection  to the set of all intersection points is an "unramified double
cover";

For every 2-cocycle $\sigma $ Donaldson constructed a so called fundamental
cycle $ D_{\sigma}$ in the space $ \cal B $$(E) $ of gauge-orbits,which is a
closed subspace of codimension 2 in $\cal B $$(E) $ (see [D 1] or the formulas
(4.18) - (4.20) in the survey article [T 1]). Then the third condition is

  3) the collection of fundamental cycles $ D_{\sigma_1},...,D_{\sigma _{d_1}}$
is in general position with respect to the strata of ends of $\cal M$$^g _1
(E)$.

  Then we can define the value of the  $Spin ^\bbbc$-polynomial as the
algebraic number of points of intersection \[ \gamma ^{g,C}_1(E)
(\sigma_1,...,\sigma_{d_l})=D_{\sigma_1}\cap...\cap D_{\sigma_{d_1}}\cap {\cal
M} ^{g,C}_1 (E) {~~~} (2.8) \]

 This definition makes sense because of
 \[ \]

{\bf Analogue of Donaldson's Lemma }. 1) If $\Sigma _1,...,\Sigma _{d_1}$ are
chosen in general position and
\[  c_2 (E) \geq \frac{3}{2} (b^+_2 + 1 ) - \frac{1}{2} c_1. C - 2\chi_C(L_0)
{~~~~~} (2.9) \]
  then the intersection  (2.8) is compact.

2) If $ g_t$ is a generall path in the space of metrics (which doesn't
intersect walls (see (2.19)), then the union of all intersections (2.8) is
smooth and compact.

3) The intersection number (2.8) depends only on the homology classes of the
$\sigma_i$'s.

 {\bf  Proof}. Since  we use the same ideas and constructions as in [D 1] we
will prove  only first statement where our constants are a little bit different
from Donaldson's. The proof of other statements is left to the reader.

  {\bf Remark}. In our applications we will use the $SO(3)$ - bundles with
$w_2(E) \neq 0$, so one doesn't even need the estimate (2.9).

 If the intersection (2.8) is not compact, then there exists a sequence of
connections  \[  \{ a_i \} \in  \bigcap_{i=1}^{d_1} D_{\sigma _i} \cap {\cal
M}^{g,C}_1
 (E) \] which after suitable gauge transformations will converge uniformly
(with bounded norm of the curvature) on $ M - \{m_1,...,m_l \} $, where
$\{m_i\}$ is a finite set of points of M, which can be regularized by an $ L^p$
- gauge transformation. The limit
connection can be regularized to an anti self dual connection

  \[  a_ {\infty } \in {\cal M}^{g,C}_1 (E')  with   \]
\[ c_1 (E') = c_1 (E),  c_2 (E')= c_2 (E) - l \ \]
 \[ \chi _C (E') =  \chi _C (E) + l \]

(because of the regularisation theorem of Uhlenbeck [ U ]).

Moreover $l \leq c_2 -\frac{1}{4} c_1^2 $ (otherwise ${\cal M}^g(E')$ is empty
by the Bogomolov inequality).

  Since every $ D_{\sigma_i}$ is closed, there are two posibilities: \[
either  \ \ a _{\infty} \in  D _{\sigma _i} {~~~~~} (2.10) \]

 \[  or {~~~}\exists {~~~}j{~~~} such{~~}  that  \ \     m_j \in  {\Sigma }_i
{~~~~~~~~} \]

Consider first the extremal case $l = c_2 - \frac{1}{4}.c_1^2 + 1$. Then the
flag diagram (2.7) gives the inequality
\[ 2 l =\#\{ m_i \in \Sigma_j \} \geq \# \{\Sigma_1,...,\Sigma_{d_1} \} {~~~}
(2.11) \]

But
\[ d = \frac{1}{2} dim_{\bbbr}{\cal M}^g(E) = 4c_2 - c_1^2 -\frac{3}{2} (b^+_2
+ 1) {~~~} (2.12)\]
\[d_1 =  \frac{1}{2} dim_{\bbbr}{\cal M}^{g,C}_1(E) = 3c_2 - 1 - \frac{1}{2}
c_1 (c_1 - C) - \frac{3}{2}(b^+_2 + 1) - 2\chi_C(L_0) \]
{}From this the inequality (2.11) is equivalent to
\[ c_2 < \frac{3}{2} (b^+_2) -\frac{1}{2} c_1.C -   2\chi_C(L_0) {~~~} (2.13)
\]
contradicting (2.9).

In the general case

                 \[a_{\infty } \in  \bigcap_{i=1}^{l'} D_{\sigma _i} {~~~~~~}
(2.14)  \]  and each surface  $ \Sigma _{l' + 1},...,\Sigma _{d_1}  $ contains
at least one point in $ \{ m_1,...,m_ l\} $.

  By the general position condition 3) \[   \bigcap_{i=1}^{l'} D_{\sigma _i}
\cap {\cal M}^{g,C}_1 (E') \neq \phi \Rightarrow   \frac{1}{2} dim_{\bbbr}
{\cal M}^{g,C}_1 (E') \geq l' {~~~~~~~}(2.15) \]

  On the other hand the flag diagram (2.7) gives the inequality \[ 2 l \geq d_1
- l' \Longleftrightarrow l'  \geq d_1 - 2 l {~~~~~~~} (2.16) \]

  From (2.15) and (2.16) we have  \[  \frac{1}{2} dim_{\bbbr} {\cal M}^{g,C}_1
(E') \geq    \frac{1}{2} dim_{\bbbr} {\cal M}^{g,C}_1 (E) - 2 l \]

that is \[ 4. c_2 - 4 l - c^2_1 -  3  (b^+_2 + 1 )  - 1 + \chi_C (E) + l \geq
 4. c_2  - c^2_1 -   3 (b^+_2 + 1 )  - 1 + \chi_C (E)) - 2 l  \]
 and this is a contradiction if l is positive.

   This proves the lemma and completes - with the usual additions (see 3.1-3.3
from [D 1]) the construction of the   Spin - polynomials
 \[ \gamma^{g,C}_1 (2,c_1,c_2) \in S^{d_1} H^2(M, \bbbz)  {~~~~~}   (2.17) \]
for a regular (in the sense the Transversality Theorem ) metric g avoiding
reducible connections.

Recall that if   $ b^+_2 = 1$,then associating to the metric g  the ray of
harmonic selfdual forms on M defines the so called period map  of the space of
Riemannian metrics to the Lobachevski space \[ K^+ \subset H^2 (M,\bbbr) /
\bbbr^+  {~~~~} (2.18) \]

The Lobachevski space $ K^+$ is divided by the collection of walls $
\{W_e = e^{\perp} \}$
\[ e \in H^2 (M,\bbbz), e \equiv c_1 mod 2, c^2_1 - 4 c_2 \leq e^2 \leq 0 {~~}
(2.19) \]
into chambers of type $(c_1,c_2)$, which form a set $\Delta$ .

Actually we can lift $ SO(3) $ - connections up to $U (2) $- connections, then
the reducibility conditions give the decomposition of our vector bundle as $
E = L_e \oplus L_{c_1 - e} $.The wall e will be important for us if
 $ \chi_C (L_e) $ or $ \chi_C (L_{c_1 - e}) $ will be positive.
Now we can compute the link of the singularity of the bordism of the moduli
spaces given by one dimensional path in the parameter space ${\cal S} \times
\Omega^1 $ and the number of points of the intersection (2.8) which disappeared
(appeared) in (from) this singularity. This number is given by the pure
topological formula actually by Porteus formula for the virtual index vector
bundle of the family of Dirac operators (the details you can see in the
forthcoming preprint of Victor Pidstrigach) .From this by the same bordism
arguments as in [D 1] and [K] we obtain the description of the dependence of
the  Spin-polynomials on the parameters:
\[   \]
 {\bf Theorem 2.1. } If $ b^+_2 = 1$  then for every pair of regular metrics
$g_1,g_2$ from the same chamber $ C\in \Delta $ \[ \gamma^{g_1,C}_1 (2,c_1,c_2)
=  \gamma^{g_2,C}_1 (2,c_1,c_2)  {~~~~~~} (2.20) \]

(if    $ b^+_2 \geq 3 $ this is true without any chamber condition).
\[  \]
  On the other hand the dependence of the Spin-polynomials on changing the  $
Spin ^ \bbbc $-structure $ C \in H_w(M) $ was given in section 1:
\[ \gamma^{g_1,C+ 2 \delta}_1 (2,c_1,c_2)  =  \gamma^{g_1,C}_1 (2,c_1- 2
\delta,c_2 -c_1. \delta +  \delta ^2 )  {~~~~}(2.21)
\] by the formula (1.12).

 \section {Algebraic surfaces }

   If M is the underlying manifold of an algebraic surface S, then there exist
the canonical  $ Spin ^ \bbbc $-structure  given by the anticanonical class $
-K_S$ (we will drop the index as long as there is no danger of confusion).

  In this case for Hodge metric  $ g_ H $ given by a polarization
 \[ H \in Pic S \subset H^2 (S,\bbbz) {~~~~~~} (3.1) \]
 the Donaldson-Uhlenbeck identification theorem gives
 \[    {\cal M}^{g_H} (E) = M^H (2,c_1,c_2 ) {~~~~~~~} (3.2) \]
where the right side is the moduli space of H - slope stable bundles on S with
Chern classes $ c_1,c_2.$

 Under the identification (3.2) $(a) = E $ we have an identification

\[ ker D^{g_H,-K}_a = H^0 (E) \oplus H^2 (E)  {~~~~~~~~} (3.3) \] \[ coker
 D^{g_H,-K}_a = H^1 (E) \] where the $H^i(E)$ denote the coherent cohomology
groups and \[ ind  D^{g_H,-K}_a  = \chi (E) {~~~~~~~~~} (3.4) \]  \[\chi _{-K}
(L_0) = \chi ( {\cal O}_S ) = p_g + 1 \]

  The subspace (2.2) then is the Brill - Noether locus  \[ {\cal M}^{g_H,-K}_1
(2,c_1,c_2) = \{ E \in  M^H (2,c_1,c_2 ) | h^1 (E) \geq -\chi (E) + 1 \} {~~~~}
(3.5) \]

  But in the  situation of surfaces the last inequality \[ h^1 (E) \geq - \chi
(E) + 1  \Longleftrightarrow  h^0 (E) + h^2 (E) \geq 1 {~~~~~~}(3.6) \]

   Hence we have a decomposition \[ { \cal M}^{g_H,-K}_1 (2,c_1,c_2) =
M_{1,0}^H (2,c_1,c_2) \cup   M_{0,1}^H (2,c_1,c_2) \]

  where the components are algebraic subvarities  \[  M_{1,0}^H (2,c_1,c_2) =
\{ E \in  M^H (2,c_1,c_2 )  | h^0 (E) \geq 1 \} {~~~~~~~~~} (3.7) \]   \[
M_{0,1} ^H(2,c_1,c_2 ) =  \{ E \in  M^H (2,c_1,c_2 )  | h^2 (E) \geq 1 \}  \]

   On the other hand the transformation  \[ E \leadsto  E^* (K) = E^* \otimes
{\cal O}_S (K)    {~~~~~~~~~~} (3.8) \]  gives the identification
\[  M^H (2,c_1,c_2 )  =   M^H (2,2 K - c_1,c_2 -c_1 .K + K^2 )
{~~~~~~~~~~~~}(3.9)  \]
 and by Serre-duality
  \[ M_ {0,1}^H (2,c_1,c_2 )  =    M_{1,0}^H (2,2 K - c_1,c_2 -c_1 .K + K^2 )
{~~~~~~~} (3.10)  \]
By this reason we have in the algebraic geometric situation two polynomials
\[  \gamma^{g_H,-K}_{1,0} (2,c_1,c_2) ,{~}  and {~}  \gamma^{g_H,-K}_{0,1}
(2,c_1,c_2)  {~~~~~~~~} (3.11)  \]
 given by the construction in section 2 with the subspaces (3.7) ( of course,
if our Hodge metric $ g_H $ is regular, then  the spaces (3.7) have the
expected dimension ).

 Now to compute the Spin-polynomial (2.17) we must sum the individual
polynomials (3.11) but here we must be careful because the natural orientations
of the components (3.7) can be different .Actually in section 5 of Ch.1 of
[P-T] following  the orientation law was proved:
\[  \]

 {\bf Orientation Rules }.1) If the number \[  1 - \chi ( E ) {~~~~} (3.12) \]
is even,then the natural  orientations of $  M_ {1,0}^H (2,c_1,c_2 )  $ and   $
 M_ {0,1}^H (2,c_1,c_2 )  $   coincide (compared with the  complex
orientation).

  2) otherwise they have different orientations.
\[   \]
  It means that
\[ c_2 = \frac{1}{2} ( c^2_1 - c_1.K) + 1  mod 2 \Rightarrow
\gamma^{g_H,-K}_1 =
 \gamma^{g_H,-K}_{1,0} +  \gamma^{g_H,-K}_{0,1}  {~~~} (3.13) \]
 and  \[  c_2  = \frac{1}{2} ( c^2_1 - c_1.K) + mod 2   \Rightarrow
\gamma^{g_H,-K}_1 =
 \gamma^{g_H,-K}_{1,0} -  \gamma^{g_H,-K}_{0,1}   \]

   On the analogy of the Non-degeneracy Theorem for the original Donaldson
polynomial we can prove
\[  \]
  {\bf Theorem 3.1}.Assume our Hodge metric $g_H$  avoids reducible connections
and \[ 1)  c_2 \geq 5 (p_g + 1) + \frac{1}{2}c_1 .K {~~~}   (see (2.9);\]

\[  2)  c_2 = \frac{1}{2} ( c^2_1 - c_1.K) + 1 mod 2  \]

\[  3)    M_ {1,0}^H (2,c_1,c_2 ) {~~~}  and {~~~}   M_ {0,1}^H (2,c_1,c_2 ) \]
 have an expected positive dimension,then \[   \gamma^{g_H,-K}_1 (2,c_1,c_2)
\neq  0{~~~~~~~~~} \]
\[    \]
 {\bf Proof.} Because of condition 1) the polynomial exists. We can choose some
smooth curve $ C $ in the complete linear system $| N H |, N \gg 0 $, such that
the restriction map \[ res_C :  M_ {1,0} ^H(2,c_1,c_2 ) \cup   M_ {0,1}^H
(2,c_1,c_2 ) \rightarrow M_C (2,c_1.C ){~~~}(3.14) \]  is an embedding (see,for
example, [T 3]).On the other hand, there is an ample divisor $ \Theta  \in Pic
 M_C (2,c_1.C ) $ and the value of Spin-polynomial on the class C is the sum of
degrees of the image subvarieties $( res_C (  M_ {1,0}^H (2,c_1,c_2 ) \cup   M_
{0,1}^H (2,c_1,c_2 ) ) $ with respect to this $\Theta $.

  It must be a sum, not a difference because of condition 2) (see (3.13)). As
in the case of the original Donaldson polynomials we are done.

  The condition 3) is very important because the properties of Spin-polynomials
depend on degrees
 \[ deg_H c_1 = H.c_1  ,{~~} deg_H K_S = H.K_S    {~~~~~~~~~~~} (3.15) \]
 with respect to the polarisation H. Namely, in contrast to the behaviour of
the Donaldson polynomials the Spin-polynomials may vanish for all values of
$c_2$:

 {\bf Lemma 3.1 }If \[ 2 K_S .H \leq  c_1.H  \leq  0 {~~~~~~~} (3.16) \]
 then
 \[  M_ {1,0} ^H(2,c_1,c_2 ) \cup   M_ {0,1}^H (2,c_1,c_2 ) = \phi \]
 and ,hence  \[ \gamma^{g_H,-K}_1 (2,c_1,c_2) = 0   \] for every $c_2$.

  {\bf Proof.} Indeed,  \[ h^0 (E) > 0 \Longleftrightarrow  \exists  s: {\cal
O}_S \rightarrow (E), s \neq 0       \] On the other hand \[  h^2 (E) > 0
\Longleftrightarrow  \exists  j: E \rightarrow    {\cal O}_S (K_S),j\neq 0 \]
but $ s \neq 0, c_1 .H \leq 0 $ contradicts the stability condition for E
and $ c_1 .H \geq 2 H.K_S $  is a contradiction to the stability condition for
E , too. We are done.

{\bf Remark.} Of course, the inequalities (3.16) are possible for rational
surfaces only. The vanising condition (3.16) is  crucial (it is actually due to
  Donaldson [D 2]).
  Of course the original Donaldson polynomials don't vanish under this
conditions as for example in the case $ S = {\bbbc}{\bbbp}^2 $ for the sequence
\[ (2,-2,c_2) , c_2 \in {\bbbz}^+ \]

  \section { Asymptotic regularity }

   Let S be a algebraic surface, H  a polarisation on S and $ c_1 \in Pic S $
a divisor class.

  {\bf Definition 4.1}A class $ c_1 \in Pic S $  is called H-semisimple, if for
any effective curve $ C \subset S $  \[ c_1.H > 2 C.H  \Longleftrightarrow
C.K_S + C^2  \leq c_1.C  {~~~~~~~} (4.1)  \]

(I would like to emphasize that the left side of the last inequality is the
degree of the canonical class on C).

   On the analogy of Donaldson's Non-degeneracy Theorem we prove
\[   \]
  {\bf Theorem 4.1.} For every H-semisimple $ c_1 \in Pic S $ with $ c_1.H > 0
$ there exists a constant $ N(H,c_1) $ such that for $c_2 \geq N (H,c_1)$
\[v.dim    M_ {1,0} ^H(2,c_1,c_2 ) = dim    M_ {1,0} ^H(2,c_1,c_2 ) > 0  {~~~}
(4.2) \]
and general point of $   M_ {1,0} ^H(2,c_1,c_2 ) $ is smooth.
 \[   \]

{\bf Proof} 7  Each $ E \in   M_ {1,0} ^H(2,c_1,c_2 )  $ has a section,that is
a non zero homomorphism \[  s: {\cal O}_S \rightarrow (E)  {~~~~}(4.3) \]

  The subscheme of zeroes of this homomorphism contains a priori subschemes of
different dimensions:

 \[ (s )_0 = C \cup  \xi {~} with {~~~} dim C = 1,   dim \xi = 0 {~~~} (4.4) \]
  Because E is H-stable , we have

 \[ 2 C.H < c_1.H  {~~~~~} ( 4.5)\]

There exists a finite set of non empty complete linear systems

\[ | 0 |, | C_1 |,...,| C_N | {~~~~~}(4.6) \]
satisfying the inequality (4.5) ( $|0|$ is the complete linear system of the
class  $0 \in Pic S $ ).

 For every $ i = 0,1,...,N $ consider the variety

 \[  GAM _ {C_i} (2,c_1,c_2 ) = \{ 0 \rightarrow {\cal O}_S (C_i) \rightarrow E
\rightarrow J_{\xi} (c_1 - C_i )  \rightarrow 0 \} / {\bbbc }^* {~~~} (4.7) \]
of all non trivial extensions up to homotheties, where $J_{\xi} $ is the ideal
sheaf of a 0-dimensional subscheme $ \xi $ (of a cluster $\xi$ for short).
(GAM  alias GAMBURGER ).

We need to prove that

 \[ dim \bigcup_{i=0}^{N}   GAM _ {C_i} (2,c_1,c_2 )  \leq v.  dim    M_ {1,0}
^H(2,c_1,c_2 )  {~~~~}(4.8) \]  and that $   M_ {1,0} ^H(2,c_1,c_2 ) \neq \phi
$.

 But under the operation $E \leadsto E(-C_i)$

\[  GAM _ {C_i} (2,c_1,c_2 )  = GAM _ 0   (2,c_1 -2 C_i,c_2 - c_1.C_i + C^2_i)
{~~~~~}(4.9) \]

 The constants  $ \{C_i^2 - c_1.C_i \} $ are bounded, hence we are done if we
prove the following
\[  \]
 {\bf Lemma 4.1}For $c_2 \gg 0$

\[ dim   GAM _ 0 (2,c_1,c_2 )  \leq   v.  dim    M_ {1,0} ^H(2,c_1,c_2 ) =   \]

\[ = 3 c_2 - 1 - \frac {c_1 (c_1 + K ) }{2} - (p_g + 1) {~~~} ( 4.10) \]
\[   \]
 (see (2.12) with $C = -K_S $ )

Note that

\[ 3 (c_2 - c_1.C_i + C^2_i ) - 1 - \frac {(c_1 - 2 C_i)( K + c_1 - 2 C_i)}{2}
- (p_g + 1) = \]

\[ = 3 c_2 - 1 -  \frac {c_1 (c_1 + K )}{2}  - (p_g + 1) + (C_i.K + C_i^2 -
c_1.C_i ) \]

and the tail is non positive due to the inequality (4.1).
\[   \]
  {\bf Proof of Lemma 4.1}The natural projection

\[ \pi :   GAM _ 0 (2,c_1,c_2 ) \rightarrow Hilb^{c_2} S {~~~~}(4.11) \]
given by sending the extension (4.7) to the cluster  $\xi $ as element of the
Hilbert scheme is surjective for big $ c_2$ . A fibre

\[ \pi ^{-1}(\xi ) = {\bbbp} Ext^1( J_{\xi} (c_1),{\cal O}_S ) =
 {\bbbp} H^1 ( J_{\xi} (c_1 + K))^*    (4.12)  \]

by Serre-duality.

 For every $ \xi \in Hilb^{c_2} S $ we have a short exact sequence

\[ 0 \rightarrow J_{\xi} (c_1 + K) \rightarrow {\cal }O_S(c_1 + K) \rightarrow
{\cal O}_{\xi} (c_1 + K) \rightarrow 0 \]
giving rise to a cohomology exact sequence

\[H^0 ( J_{\xi} (c_1 + K) )\rightarrow H^0({\cal }O_S (c_1 + K)) \rightarrow
{\bbbc}^{c_2} \rightarrow H^1 ( J_{\xi} (c_1 + K) )\rightarrow H^1({\cal }O_S
(c_1 + K)) \rightarrow 0  {~}(4.13) \]

 Moreover,

\[  h^0 ( J_{\xi} (c_1 + K) ) = 0 \Rightarrow dim  {\bbbp} H^1 ( J_{\xi} (c_1 +
K)) = c_2 - \chi  (   {\cal }O_S(c_1 + K) ) - 1 = \]

\[ = c_2 - 1   -  \frac {c_1 (c_1 + K )}{2}  - (p_g + 1)  {~~~}(4.14) \]

  Consider the subvariety

\[ \Delta = \{\xi \in Hilb^{c_2} S | h^0 ( J_{\xi} (c_1 + K) ) > 0 \} {~~~}
(4.15) \]
 It is easy to see that

\[ dim \Delta \leq c_2 + dim | c_1 + K | = c_2 +  h^0({\cal }O_S(c_1 + K) ) - 1
\]
On the other hand from (4.13) we have

\[    dim    \pi^{-1} (\xi ) =  {\bbbp} H^1 ( J_{\xi} (c_1 + K)) \leq
h^1({\cal }O_S(c_1 + K) ) +c _2 \]
Hence

\[    dim     \pi^{-1} (\Delta )  \leq 2 c_2 - 2 +     h^0({\cal O} _S(c_1 + K)
) +  h^1 ({\cal }O_S(c_1 + K) ) \]

and

\[ c_2 >  2  h^0({\cal }O_S(c_1 + K) ) \Rightarrow   dim  \pi^{-1} (\Delta ) <
3 c_2 - 2 -  \frac {c_1 (c_1 + K )}{2}  - (p_g + 1) . \]

 This proves  Lemma 4.1.To finish the proof of Theorem 4.1 we prove

 {\bf Lemma 4.2 }.If $ c_1.H >0 $,then for $c_2 \gg 0$ \[   M_ {1,0}
^H(2,c_1,c_2 )  \neq \phi \]

and hence by Theorem 4.1 it has the expected dimension.

 {\bf Proof }.We need to prove that for generic $ \xi \in Hilb^{c_2} $ and $
c_2 \gg 0 $ any non trivial extension
\[   0 \rightarrow {\cal O}_S  \rightarrow E \rightarrow J_{\xi} (c_1  )
\rightarrow 0  \]
is  H-stable.
Twisting E by $ (-c_1) $ we have
\[   0 \rightarrow {\cal O}_S (-c_1)  \rightarrow E (-c_1)\rightarrow J_{\xi}
  \rightarrow 0  \]

The hypothetical destabilizing line bundle must be of type $ {\cal }O_S(-
C)$,where C is an effective curve subject to the inequality (4.1) and the
cluster  $\xi $  must be supported on this effective curve.

But the collection (4.6) of such curves is finite and for  $ c_2 \gg 0 $ (as in
Theorem 4.1) a generic $ \xi $ is not contained in any curve in this collection
of complete linear systems (see (4.15)).

 {\bf Definition 4.2 }. A class $c_1 \in Pic  S $ is called H-simple, if it is
semisimple and the class $ 2K_S - c_1 $ is semisimple too.

 As a corollary of Theorems 3.1 and 4.1 we provide

 {\bf Theorem 4.2}.Assume that our Hodge metric g avoids reducible connections
. Then for every H-simple  $c_1 \in Pic  S $ with $ c_1.H >0 $
there exists a constant  $ N(H,c_1)$ such that  for \[ c_2 \geq N(H,c_1) , c_2
= \frac{1}{2} ( c^2_1 - c_1.K) + p_g mod 2  {~~} (4.16)\]
\[  \gamma^{g_H,-K}_1 (2,c_1,c_2) \neq 0   \]

 At last (but not  at least ) we need to explain what we have to do if our
Hodge metric does not avoid reducible connections.Certainly in case when $rk
Pic S > 1 $ we may use the following extremely useful trick:

 {\bf Definition 4.3}.A polarization $ H^{\varepsilon} $ is called  close to H
if the ray $ {\bbbr}^+.H^{\varepsilon} $ in $K^+ $ (2.18) is close to the ray
 $ {\bbbr}^+.H $ in Lobachevski metric.

 {\bf Lemma 4.3}. If a class $ c_1 \in Pic S $ is a H-simple ,then for a
polarisation $ H^{\varepsilon} $ sufficiently close to H
\[   1) {~~~} c_1 {~~~} is {~~~}  H^{\varepsilon } -simple {~}too.\]

\[ 2) {~}  2 K_S .H <  c_1.H  <  0 \Longrightarrow   2 K_S .H ^{\varepsilon}<
c_1.H ^{\varepsilon}  <  0 \ {~~~~} (4.17) \]

  3) for every polarisation  H there exists a sufficiently close to H
polarisation $ H^{\varepsilon}$  such that the Hodge metric $
g_{H^{\varepsilon}} $ avoids the reducible connections.

 {\bf Proof }.For a sufficiently close polarisation  $ H^{\varepsilon} $

\[ 2 C.H < c_1. H \Longrightarrow  2 C. H^{\varepsilon} < c_1. H^{\varepsilon}
{~~~~~} (4.18) \]

Hence the collection of linear systems (4.6) for   $ H^{\varepsilon} $
is the same as for H and we have 1) and 2).

To prove 3) it is enough to remark that the set of  rays of polarisations is
dense in the projectivisation of the K\"ahler cone in $ K^+ $ and the set of
walls is discrete and locally finite.

  In the last section we consider three very simple examples to show how we can
use Theorem 4.2 and Lemma 3.1 to distinguish the underlying smooth structures
of rational surfaces and surfaces of general type.

 \section { Applications }.

  For the beginning we prove

   {\bf Theorem 5.1 }.If an algebraic surface S is diffeomorphic to   $
{\bbbc}{\bbbp}^2 $ then   $ S = {\bbbc}{\bbbp}^2 $ (as algebraic surface).

 {\bf Proof}.Let

\[   f:  {\bbbc}{\bbbp}^2  \rightarrow  S {~~~~~~~~~~~} (5.1) \]
be a diffeomorphism, $ h \in Pic S $ be the positive generator of $Pic S$ ($h^2
=1$).Then (using the real anti involution of $  {\bbbc}{\bbbp}^2 $ if
necessary) we may consider the case when

\[          f^*  \  (h) = l     \]

is the class of the line on  $ {\bbbc}{\bbbp}^2 $.

For the canonical class we have $K_S = 3 h$ (otherwize  $K_S =-3h $ and S is
rational).

Then the increment of the canonical class with respect to f is

\[  \delta _f (K) =  \frac {  f^*(K_S) -K_{   {\bbbc}{\bbbp}^2  } }{2} = -K_{
{\bbbc}{\bbbp} }^2  {~~~~} (5.2) \]

Consider the following topological type of vector bundles on S
\[    (2,h,c_2),  c_2 \gg 0  {~~~~~~~~~~} (5.3)  \]

  {\bf Lemma 5.1 }.For all $c_2$ the Spin-polynomial

\[   \gamma^{g_h,-K}_1 (2,h,c_2) = 0 .     \]

  {\bf Proof}.The operation $f^*$ gives the equality

\[  \gamma^{g_h,-K}_1 (2,h,c_2) =  \gamma^{f^*(g_h),K_{\bbbc \bbbp ^2}}_1
(2,l,c_2) \]

By the equality (2.20)

\[   \gamma^{f^*(g_h),K_{\bbbc \bbbp ^2}}_1 (2,l,c_2) =
\gamma^{g_{F-S},K_{\bbbc \bbbp ^2}}_1 (2,l,c_2)   \]

where $g _{F-S} $ is the Fubini-Study metric on $ \bbbc \bbbp ^2 $.

By the equality (2.21)

\[   \gamma^{g_{F-S},K_{\bbbc \bbbp ^2}}_1 (2,l,c_2) =
\gamma^{g_{F-S},-K_{\bbbc \bbbp ^2}}_1 (2,2 K_{\bbbc \bbbp^2}+l,c_2+ 6)=
\gamma^{g_{F-S},-K_{\bbbc \bbbp ^2}}_1 (2,- 5l,c_2 + 6) \]

The first Chern class $c_1 = -5l $  satisfies the inequality (3.16) and by
Lemma 3.1 we are done.

To provide a contradiction to the existence of an f (5.1) we prove

 {\bf Lemma 5.1' }.On S the class h is h-simple. Hence by Theorem 4.2 if $
c_2$ is odd then

\[   \gamma^{g_h,-K}_1 (2,h,c_2) \neq  0            \]

  {\bf Proof} For h the set of linear systems (4.6) is $|0|$.Hence h is
h-semisimple.

  For the class $2K_S - h = 5h $ the set (4.6) is

\[  |0|, |h|, |2h|   \]

For these classes we have

\[ h^2 + h.K_S =4 <5 , (2h)^2 + 2h.K_S = 10 \leq 10  \]
and thus $2K_S -h$ is semisimple too . This implies that h is simple.We are
done by Theorem 4.2.

  Let $ F_1$ be the projective plane blown up in one point (Hirzebruch surface
of number 1).

  {\bf Theorem 5.2}.If an algebraic surface S is diffeomorphic to $F_1$ then $
S = F_{2n+1} $ that is the odd Hirzebruch surface.

 {\bf Proof}.Let

\[   f: F_1 \rightarrow S  {~~~~~~} (5.4) \]

be a diffeomorphism.We can find a basis h,e in $Pic S$ such that

\[ h^2 = 1, e^2 = -1, h.e = 0, K_S = 3h - e   {~~~~} (5.5)\]

Again it is sufficient to consider the case when

\[ f^* (h) = l, f^* (e) = E {~~~~~~~~~} (5.6) \]

where l is the class of line and E is the exceptional divisor on $F_1$.Hence

\[  K_{F_1} = -3l + E = f^*(-K_S)  {~~~~~} (5.7) \]

Then the increment of the canonical class

\[ \delta _f K = -K_{F_1}        \]

We only need to investigate the case when S is a surface of the general type
and minimal.

Then $ K_S $ is a polarisation on S and $f^* (K_S)=-K_{F_1}$ is a polarisation
on $F_1$.
Let $ H $ be a polarisation on S sufficiently close to $K_S$ such that

\[ f^*(H) = H_1  \]

is a polarisation on $F_1$ sufficiently close to $(-K_{F_1})$.

  {\bf Lemma 5.2}.For all $c_2$ the Spin-polynomial

\[   \gamma^{g_H,-K_S} _1 (2,h,c_2) = 0 .     \]

{\bf Proof}.The operation $f^*$ gives the equality

\[  \gamma^{g_H,-K_S}_1 (2,h,c_2) =  \gamma^{f^*(g_{H_1}),K_{F_1}}(2,l,c_2) .
\]

But the metrics   $f^*(g_H) $ and $  g_{ H_1} $  on $F_1$ are contained in the
same chamber. More precisely they have the same image of the period map .Then
by (2.20)

\[  \gamma^{f^*(g_H),K_{F_1}}(2,l,c_2) = \gamma^{g_{H_1},K_{F_1}}(2,l,c_2)   \]

Moreover

\[ \gamma^{g_{H_1},K_{F_1}}(2,l,c_2) =
\gamma^{g_{H_1},-K_{F_1}}(2,-5h+2l,c_2+5) \]

The first Chern class $c_1 = -5h +2l$ satisfies the inequality (3.16). Thus by
Lemma3.1 we are done.

To provide a contradiction to  existence of f (5.4) we prove

 {\bf Lemma 5.2'}.On S the class h is $K_S$-simple.

{\bf Proof} Let C be an effective curve C of the form $ C=xh -ye $. We will
check whether the inequality (4.1) does hold for C : First

\[2 \leq 2.C.K_S \leq h.K_S =3 \Rightarrow C.K_S = 1 \Rightarrow y=3x-1  \]

Then

\[C^2 = x^2 -y^2 =-8x^2 + 6x -1,C^2 + C.K_S = -8x^2 +6x <C.h = x \]

for every $x \in \bbbz $.Hence  h is $K_S$-semisimple.

Now for $c_1 = 2 K_S - h = 5h - 2l $ we have

\[ 2 \leq 2C.K_S \leq 13  \Rightarrow C.K \in \{1,...,6\} \Rightarrow y = 3x -
\{ 1,...,6\}  \]

Then

\[ C^2 + C.K_S = -8x^2 + 6x .\{1,...,6\} - \{0,2,6,12,20,30\} \]

\[ C. (5h-2e) = - x + 2.\{1,...,6\} \]

It is easy to check that for any of the six cases the inequality (4.1) holds.
Hence $5h - 2e $ is $K_S$-semisimple ,too, and we are done.

  At last let us go to the Hirzebruch Problem.

{\bf Theorem 5.3}.On $ S^2 \times S^2 $ there exists the unique algebraic
structure $ Q = \bbbc \bbbp^1 \times \bbbc \bbbp^1 $ up to the elementary
transformations to the even Hirzebruch surface $F_{2n}$.

{\bf Proof}.In this case for any topological type $(2,c_1,c_2) $ of vector
bundle the virtual dimension  of $ M^H(2,c_1,c_2) $ is odd.We will use a simple
trick:

Let

\[    f:Q \rightarrow S  {~~~~} (5.8) \]

be a diffeomorphism. We can find basis $ h_+,h_- $ in $Pic S$ such that

\[ h_+^2 = h_-^2 = 0, h_+.h_- = 1,K_S = 2h_+ + 2h_-  {~~~} (5.9) \]

such that $f^* (h_+)= h'_+$ , $f^*(h_-)=h'_- $ is the standard basis of $ Pic
            \bbbc \bbbp^1 \times \bbbc \bbbp^1 $ and $f^*(K_S) = -K_Q $.

We only need to consider the case when $K_S$ is nef.

Let's blow up a point p on Q and $ f(p) $ on S. Then the diffeomorphism f (5.8)
can be extended to a diffeomorphism
\[ \tilde{f} : \tilde{Q} \rightarrow  \tilde {S} {~~~~} (5.8') \]
and
\[ K_{\tilde{Q}} = K_Q + E' ,  K_{\tilde{S}} = K_S + E  {~~~} (5.10) \]
where E and E' are the respective exceptional curves.
The increment of the canonical class with respect to $\tilde{f}$ is
\[ \delta _{\tilde{f}} K = -K_Q     \]
The divisor class $ (h_+ + h_-) \in Pic \tilde{S} $ is nef and we consider a
polarisation H on $\tilde{S} $ sufficiently close to $  (h_+ + h_-) $ such that
 $ \tilde{f}^*(H) = H'$ is a polarisation on $\tilde{Q} $ sufficiently close to
  $ (h'_+ + h'_-) $ .

 {\bf Lemma 5.3}.For all $c_2$  the Spin-polynomial
\[\gamma ^{g_H,-K_S -E}_1 (2,h_+ +h_- +E,c_2) = 0  \]

 {\bf Proof} As usual by (2.20) and (2.21)
\[\gamma ^{g_H,-K_S - E}_1 (2,h_+ +h_- +E,c_2) = \gamma ^{g_H',-K_Q - E'}_1
(2,-3 (h'_+ + h'_- + -E') ,c_2 + 2)  \]
The first Chern class $ c_1 = -3 (h'_+ + h'_- -E') $ satisfies the inequality
(3.16) so by Lemma 3.1 we are done.

 To provide a contradiction to the existence of $\tilde{f} $ (and hence of f
(5.8)) as before we prove

{\bf Lemma 5.3'}. On $\tilde{S} $ the divisor class
$  h_+ +h_- +E $ is H - simple.

{\bf Proof}.For an effective curve C
\[ 0 \leq C (h_+ +h_-) <1 \Rightarrow C = m(h_+ -h_-) + nE  \]
But then
\[ C.K_{\tilde{S}} + C^2 = -2m^2 - n^2 -n \leq C.(h_+ +h_- +E) = -n \]
and $  h_+ +h_- +E $ is H - semisimple.

Now for $  2 K_{\tilde{S}} - c_1 =3 (h_+ + h_-)+E $ let

\[ C = xh_+ +  m(h_+ -h_-) + nE . \]

Then
\[0 \leq C (h_+ +h_-) <3 \Rightarrow  x =\{0,1,2\} \]
and
\[  C.K_{\tilde{S}} + C^2 =  \{0,2,4\}m - 2m^2 - n^2 +\{0,2,4\} -n ,{~}
C.(3h_+ + 3h_- +E) = \{0,3,6\} -n \]
 It is easy to check that the right side of the inequality (4.1) holds for all
n and m .Hence $h_+ +h_- +E$ is H-simple.
The reader may continue these purely arithmetical investigations himself. Good
luck!.
\[ \]

  {\bf Reference }

 [D 1] S.Donaldson  "Polynomial invariants for smooth 4-manifolds" ,Topology 29
	          (1990),257-315. \\

 [D 2] S.Donaldson  "Differential topology and complex varieties",Arbeitstagung
 Proceedings 1990,MPI (1990). \\

 [F-U] D.Freed,K.Uhlenbeck "Instantons and four-manifolds" M.S.R.I.Publ.
Springer, New York, 1984. \\

 [K] D.Kotschick    "SO(3)-invariants for 4-manifolds with $b^+_2 =1$",
Proc.London Math.Soc.3:63(1991),426-448. \\

 [P-T] V.Pidstrigach A.Tyurin"Invariants of the smooth structures of an \ \ \
      algebraic surface arising from Dirac operator."Izv.AN SSSR
Ser.Math.,56:2(1992),279-371.(english translation:Warwick preprint 22 (1992).
\\

 [T 1] A.Tyurin     "Algebraic geometric aspects of smooth structure.1.The
Donaldson polynomials." Russian Math.Surveys 44:3 (1989),113-178. \\

 [T 2] A.Tyurin     "A slight generalization of the theorem of
Mehta-Ramanathan."Springer  LNM 1479 (1989),258-272. \\

 [T-3] A.Tyurin     " Spin-polynomial invariants of smooth structures of
algebraic surfaces."forthcoming Izv.AN Russia Ser.Math.,57:2(1993) \\

 [U]  K.Uhlenbeck   " Connections with $L^p$ bounds on curvature" Comm.Math .
Phys.83 (1982),31-42.

\end{document}